\magnification =1098

\hsize 16truecm
\vsize 23truecm
\hfuzz =11pt
\scrollmode

\def \a{{\alpha}}
\def \b{{\beta}}

\def \e{{\varepsilon}}

\def \g{{\gamma}}
\def \G{{\Gamma}}

\def \o{{\omega}}
\def \O{{\Omega}}
\def \p{{\varphi}}
\def \t{{\vartheta}}
\def \m{{\mu}}
\def \s{{\sigma}}

\def \aa{{d}}
\def \A{{\cal A}}

\def \E{{\bf E}\, }
\def \F{{\cal F}}

\def \NNN{{\cal N}}

\def \P{{\bf P}}
\def \q{{\quad}}

\def \qq{{\qquad}}
\def \R{{\bf R}}

\def \T{{\bf T}}
\def \ua{{\underline{a}}}

\def \ud{{\underline{d}}}

\def \uz{{\underline{z}}}

\def \Z{{\bf Z}}
\def \zz{{\cal Z}}

\def \noi
{{\noindent}}

\def\qed{\hbox{\vrule height 6pt depth 0pt width
6pt}}
\def\cqfd{\hfill\penalty 500\kern 10pt\qed\medbreak}

\font\gem=
cmcsc10 at 10 pt
\font\gum= cmbx10 at 11 pt

scaled\magstep3
\font\ph=cmcsc10  at  10 pt
\font\iit =cmmi7 at
7pt

\font\ph=cmcsc10  at  10 pt
\font\phh=cmcsc10
at  8 pt
\font\ggum= cmbx10 at 11 pt
\font\gum= cmcsc10 at 11 pt
\font\gem=
cmcsc10 at 10 pt
\font\gum= cmbx10 at 12 pt

\font\ggum= cmcsc10 at 11 pt
\font\ggum= cmbx10 at 11 pt

\def\ddate
{\ifcase\month\or January\or
February\or March\or April\or May\or June\or
July\or
August\or September\or October\or November\or December\fi\
{\the\day},\
{\sevenrm\the\year}}


\centerline{{\ggum   On
the Supremum of Random Dirichlet Polynomials}}
\vskip 1 cm

 \centerline{Mikhail {\ph Lifshits}\ and\
Michel {\ph Weber}}
 \vskip 0,6cm

{\leftskip =1,9cm \rightskip=1,9cm

\noindent {\bf Abstract: \it    
We study the supremum  of some random
Dirichlet polynomials 
$D_N(t)=\sum_{n=2}^N\e_n \aa_n n^{-\s - it}$,
where $(\e_n)$ is a sequence of independent Rademacher random variables,
the weights $(d_n)$ are multiplicative and $0\le \s <1/2$.  
The particular attention is given to the polynomials 
$\sum_{n\in {\cal E}_\tau}\e_n n^{-\s - it}$, 
${\cal E}_\tau=\{2\le n\le N\!: \! P^+(n)\le
p_\tau\}$, $P^+(n)$ being the largest prime divisor
of $n$. We obtain sharp upper and lower bounds for supremum expectation
that extend the optimal estimate of Hal\'asz-Queff\'elec
$$ \E\, \sup_{t \in \R} \left|\sum_{n=2}^N
\e_n n^{-\s - it}\right| \approx 
{N^{1-\s}\over \log N}.
$$
Our approach in proving these results is entirely based on 
methods of stochastic processes, in particular the metric
entropy method.}
\par}\footnote
{}{{\iit Date}\sevenrm : on
\ddate\par
\vskip -2pt {\iit AMS}
\ {\iit
Subject}\ {\iit Classification} 2000:  Primary 30B50, 26D05, Secondary
60G17.\par \vskip -2pt
  {\iit Keywords}:  Dirichlet polynomials,
Rademacher random variables, metric entropy method. \par }
\vskip
0,6cm
\rm

\noindent{\gum 1. Introduction and main results}
\medskip
\noindent

\noi Let $
 \{\aa_n, n\ge 1\}$ be a sequence of real
numbers.
Let $s=\s + it$ denote a complex number.
The study of the supremum
of the Dirichlet polynomials
$$ P(s)= \sum_{n=2}^N  \aa_n n^{-s}
$$
over
lines $\{s=\s+it, \ t\in \R\}$ is naturally related to that 
of
corresponding Dirichlet series, via the abscissa of uniform convergence
$$
\s_u=\inf\Big\{\s : \sum_{n=2}^\infty  \aa_n n^{-\s - it}\
\hbox{converges
uniformly over $t\in\R$} \Big\},
$$
through the relation
$$\s_u=
\limsup_{N\to \infty} {{\log  \,  \displaystyle
\sup_{t\in \R}
\big|\sum_{n=2}^N  \aa_n n^{-it}\big|} \over \log N}\  .
$$
One can refer
to Bohr [B], Bohnenblust and Hille [BH], Helson [H], 
Hardy and Riesz [HR],
Queff\'elec [Q3] for this background
and related results. This of course,
basically justifies the
 investigation of the supremum of Dirichlet
polynomials
(see for instance Konyagin and Queff\'elec [KQ]).

The following classical reduction step  enables to replace the
Dirichlet polynomial by some relevant {\it trigonometric}
polynomial. In order to recall this reduction, we introduce the
necessary notation. Let
$2=p_1<p_2<\ldots$ be the sequence of all
primes. If $  n=\prod_{j=1}^\tau
p_j^{a_j(n)}$, we write
$\ua(n)=\big\{a_j(n), 1\le j\le \tau\big\}$. Let
$\pi(N)$ denote, as usual, the number of prime numbers 
that are less or equal to $N$. 
Finally, let $\T=[0,1[=\R/\Z$ be the torus. Let us fix $N$,
put $\m =\pi(N)$, and define, for 
$\uz= (z_1,\ldots, z_\m) \in \T^\m$,
$$ Q(\uz)= \sum_{n=2}^N
 \aa_n n^{-\s}e^{2i\pi\langle
\ua(n),\uz\rangle} ,
$$
The famous H. Bohr's observation ([Q1-3]),
states
that
 $$\sup_{t\in \R} \big|P(\s +it)\big| =\sup_{\uz \in
\T^\m}
\big|Q(\uz)\big|\ .
\eqno(1.1)$$
This indeed follows
straightforwardly from Kronecker's Theorem
(see [HW], Theorem 442,
p.382).

A parallel study is also developed for
{\it random} Dirichlet polynomials and  {\it random}
Dirichlet series in the papers of Hal\'asz [Ha1-2],
Queff\'elec [Q1-3], Bayart, Konyagin and Queff\'elec [BKQ],
Kahane [K] and of Yu [Y1-3],[STY], Hedenmalm and Saksman [HS].   
Such investigations concerning random  Dirichlet series (as well as  random
power  series) go back to earlier works of Hartman [Har], Clarke [C],
Dvoretzky and Erd\" os [DE1-2], Dvoretzky and Chojnacki [DC].

%
%
%
%
 
Let $\e= \{\e_i, i\ge 1\}$  be (here and throughout the whole
paper) a sequence of independent Rademacher random variables
($\P\{ \e_i=\pm 1\} =1/2$) defined on a basic
probability space $(\O, \A, \P)$.
%
%
%
%
%
 
%
%
Consider the random Dirichlet polynomials
$$  D(s) =\sum_{n=2}^N\e_n \aa_n n^{-\s - it}.
\eqno (1.2)$$

\noi  When $\aa_n\equiv 1$, some results about the suprema are known.
If $\s=0$, then for some absolute
constant $C$, and all integers $N\ge 2$
$$  C ^{-1} {N  \over   \log  N
}\le
\E\, \sup_{t \in \R}|\sum_{n=2}^N \e_n  n^{  - it}|
\le  C {N\over
\log N}\ .
\eqno (1.3)$$
This has been proved by Hal\'asz (see [Q2-3]). In
[Q2-3]
(see also [Q1] for a first result), Queff\'elec extended Hal\'asz's
result
to the range of values
$0\le \s<1/2$; and provided a probabilistic
proof of the original one, using Bernstein's inequality for polynomials,
properties of complex Gaussian processes and the sieve method introduced by
Hal\'asz. He obtained that for some constant $C_\s$ depending on $\s$
only, and all integers $N\ge 2$
$$  C_\s ^{-1} {N^{1-\s}  \over   \log  N }\le
\E\,
\sup_{t \in \R}|\sum_{n=2}^N \e_n  n^{-\s  - it}|  \le
C_\s
{N^{1-\s} \over \log  N }\ .\eqno (1.4)
 $$
This   in fact admits a stronger form
$$  C_\s ^{-1}  \le \E\, \sup_{N\ge 2}
\sup_{t \in \R}
{|\sum_{n=2}^N \e_n  n^{-\s  - it}| \over  N^{1-\s}     (\log  N)^{-1}   }
\le    C_\s .
\eqno (1.4')
$$ 
A proof is given at the end of Section 4. We shall hereafter use 
and simplify Queff\'elec's probabilistic argument,
notably reducing the proof of the   upper bound part to the
study of suitable real Gaussian processes (which can be easily reduced
to a single one). Further, we will not use Bernstein's inequality, unlike in both
previous proofs. A simple metric entropy argument is indeed sufficient,
making the proof entirely based upon stochastic processes methods.

By developing this approach, we will also study the case when
the $\aa_n$'s are not constant and random Dirichlet polynomials are
supported by other sets than intervals of integers $[2,N]$.
At this regard, we consider the following natural extension.
For any integer $n>1$, let $P^+(n)$ denote the largest prime divisor
of $n$. Let 1$\le M<N$ be two positive integers and
define
$$S(N,M) = \big\{ 2\le n\le N : P^+(n)\le M\big\}. $$
Since
$S(N,N)=[2,N] $, these sets  naturally generalize
the notion of interval of integers.   By using the standard notation
$$\Psi(N,M):= \sharp(S(N,M)),
$$
$u=(\log N)/\log M$,  we have
([T], Theorem 6 p.405)
 $$ \Psi^*(N,M):= {\Psi(N,M)\over N} =\rho(u)
  +
{\cal O} \big( {1\over \log  y} \big) , \eqno (1.5)
$$
uniformly for $x\ge y\ge 2$,  where   $\rho(u)$ is the Dickman
function, namely the unique continuous function on $[0,\infty[$,
having a derivative on $]0,\infty[$, and such that
$$\line{$ \qq\qq\qq\qq\qq\qq \left\{  \matrix{
\rho(v)&\ \ \
=1,\qq &(0\le v\le 1),\cr
 v\rho'(v)+ \rho(v-1)&\ \ \ =0, \qq  &(v>1). \cr
}  \right.  \hfill$}
$$  
It is known that $\rho(u)$ is a decreasing positive function and that
$\log\rho(u)\sim - u\log u$, as $u\to\infty$. In other words, $\rho$ decreases
as fast as the inverse of Gamma function.  By setting $M=N^\e$ in $(1.5)$  we see that
$\Psi(N,N^\e )\sim N\rho(\e^{-1})$
for any fixed $0<\e\le 1$.

In view of (1.5), we sometimes refer to
$\Psi^*$ as to Dickman-type function.
\smallskip

Fix some positive
integer $\tau\le \pi(N)$, and recall that $p_1<p_2<\ldots$
is the sequence of primes. Put
$$
{\cal E}_\tau= {\cal E}_\tau(N)=
\big\{ 2\le
n\le N : P^+(n)\le p_\tau\big\}.
$$
Note that for $\m=\pi(N)$ we have
${\cal E}_\m=\{2,\ldots,N\}$.

The ${\cal E}_\tau$-based Dirichlet polynomials were already 
considered in [Q3]. One motivation for considering them, 
related to Rudin-Shapiro problem, will be explained later.
\medskip

We begin with a result that contains both above
mentioned estimates (1.3) and (1.4). 
\medskip

\noi {\gem Theorem 1.1.}
{\it
a) Upper bound. Let $0\le\s < 1/2$. Then there exists a
constant $C_{\s}$ such that for any integer $N\ge 2$
it is true that
 $$
\E\, \sup_{t \in \R} \big|
   \sum_{n\in {\cal E}_\tau}\e_n  n^{ -\s -
it}\big|\le
\cases{
 C_{\s} {N^{1/2-\sigma}\tau^{1/2}\over (\log N)^{1/2}}
&,\ if \ $N^{1/2} \le \tau\le N,$  \cr
& \cr
 C_{\s} {N^{3/4-\sigma}
\over (\log N)^{1/2}}
   &,\ if \  ${ N^{1/2}\over \log N}\le
\tau \le
               N^{1/2},$ \cr
& \cr
 C_{\s}  
 \
N^{1/2-\sigma}\tau^{1/2}
   &,\ if \  $1\le \tau \le
{N^{1/2} \over \log N}\ .$
} 
 $$

b) Lower bound. Let $0\le\s < 1/2$. Then
there exists a constant
$C_\s$ such that for every $N\ge 2$,
$$
\E\,
\sup_{t \in \R} \big|
   \sum_{n\in {\cal E}_\tau}\e_n  n^{ -\s -
it}\big|
\ge {C_\s\ N^{1/2-\s} \tau^{1/2} \over (\log\tau)^{1/2}}\ \cdot
\Psi^*\Big({N \over p_\tau}\, , p_{\tau/2}\Big)^{1/2} .
$$
}
\bigskip 

{
\bf Sharpness of the result}.\  It is instructive to compare the lower and
upper bounds obtained in Theorem 1.1.

\noi Consider three cases, as in the
upper bound of this theorem:
\medskip

Case I: \ $N^{1/2} \le \tau\le
N.$

Here the Dickman function vanishes from the lower bound and we
have
$\log \tau\sim \log N$. It follows from the theorem
$$
 C_{1}(\s)
{N^{1/2-\sigma}\tau^{1/2}\over (\log N)^{1/2}}
\le
E\, \sup_{t \in \R}
\big|
   \sum_{n\in {\cal E}_\tau}\e_n  n^{ -\s - it}\big|
   \le
C_{2}(\s) {N^{1/2-\sigma}\tau^{1/2}\over (\log N)^{1/2}}  
\ .
$$
Thus our bounds are optimal.
\medskip

Case II:\ ${N^{1/2}\over \log N}\le \tau \le
N^{1/2}$.

Again the Dickman function vanishes from the lower bound
and we have $\log \tau\sim \log N$. Thus
$$
C_1(\sigma)\
{N^{1/2-\s}
\tau^{1/2}
\over (\log N)^{1/2}}\
\le
 \E\, \sup_{t \in \R} \big|
\sum_{n\in {\cal E}_\tau}\e_n  n^{ -\s - it}\big|
  \le C_2(\s)
{N^{3/4-\sigma}
         \over (\log N)^{1/2}}\ .
$$
The ratio of the right and the left hand side satisfies
$$ 1\le \ {N^{1/4}\over \tau^{1/2}}\  \le
(\log N)^{1/2}\ .
$$
Thus a logarithmic gap appears.
\medskip

Case III:
$1\le \tau \le {N^{1/2} \over \log N}$\ .

Assume first that $\tau\ge N^{\e}$ for some fixed $\e>0$,
necessarily with $\e<1/2$. Then the Dickman
function produces in the lower bound just an extra constant depending on
$\e$. 
We have
$$
C_1(\s,\e)\ {N^{1/2-\s} \tau^{1/2}
\over
(\log\tau)^{1/2}}\
\le
\E\, \sup_{t \in \R} \big|
   \sum_{n\in {\cal
E}_\tau}\e_n  n^{ -\s - it}\big|
\le
   C_{2}(\s)  \
N^{1/2-\sigma}\tau^{1/2}\ .
$$
The gap is still of the logarithmic
order:
$$ 1\le (\log\tau)^{1/2}\  \le (\log N)^{1/2}\ .
$$

One should notice that an upper etimate 
$C \, N^{1/2-\s} (\tau \,\log\log N)^{1/2}$
slightly weaker than our bound in Case III was obtained in [Q3].

It is also worth of mentioning that our approach to the lower bounds
is very different from that in the preceding works [Q3], [KQ]
based on {\it deterministic} estimates valid for any polynomial,
see e.g. lower bound in (1.6) below. It would be interesting to check
whether the optimisation of parameters in deterministic estimates
enables to this approach to compete with our lower bound on 
the whole range of $\tau$.
\bigskip

Unfortunately, if $\tau$ is relatively small, namely
$\log\tau\ll\log N$, the gap between the upper and the lower bounds
in Theorem 1.1 becomes rather significant due to the small
factor $\Psi^*$
in the lower bound. Our next result, although being
not optimal, shows that the presence of $\Psi^*$ is really crucial.
\medskip

\noi {\gem Theorem 1.2.}
{\it Let $0\le\s < 1/2$. Then there exists a constant $C_{\s}$ such
that for any integer $N\ge 2$ and $\tau>\exp\{(\log\log N)^2\}$
it is true that}
$$ 
{ N^{1/2-\s} \tau^{1/2} \Psi^*\Big({N \over p_\tau} ,
p_{\tau/2}\Big)^{1/2} 
\over C_\s (\log\tau)^{1/2}}\ 
\le
\E\, \sup_{t \in
\R} \big|
   \sum_{n\in {\cal E}_\tau}\e_n  n^{ -\s - it}\big|
\le 
C_{\s}
N^{1/2-\sigma}\tau^{1/2}\Psi^*\Big({N \over p_\tau^2},
p_{\tau}\Big)^{1/2}
\ . 
$$
\bigskip

{\bf Estimates of $\ell_1$-type.}
\
The reader familiar with evaluation of Rademacher processes may wonder whether 
the brutal $\ell_1$-estimates 
$$
\E\, \sup_{t \in \R} \big| \sum_{n\in {\cal E}_\tau}\e_n  n^{ -\s - it}\big|
\le  \sum_{n\in {\cal E}_\tau}  n^{ -\s} := L(N,\tau)
$$
are useful at least in some zone of parameters. In our context the answer is negative.
Actually, one can show that
$$
L(N,\tau)\ge c\ N^{1-\s}\Psi^*\big(N,p_\tau \big) \sim c \ N^{1-\s}\rho\Big({\log N\over \log p_\tau}\Big).
$$
This is too much for good upper bounds, as one can see from two following examples. The first one handles
large $\tau$ and the second one deals with small $\tau$.

1) Let $\p_\tau\sim N^h$ with $1/2<h\le 1$. Then we see that
$$  
L(N,\tau)\ge c(h)\  N^{1-\s}
$$
while the upper bound from Theorem 1.1 yields a better estimate
$$
\E\, \sup_{t \in \R} \big|\sum_{n\in {\cal E}_\tau}\e_n  n^{ -\s -it}\big|\le
C_{\s} 
{N^{1/2-\sigma}\tau^{1/2}\over (\log N)^{1/2}}
\approx
C_{\s} 
{N^{(1+h)/2-\sigma}\over \log N}.
$$
The gap between the two upper bounds is at least  logarithmic for $h=1$ and polynomial for $h<1$.

2) Let $\p_\tau\sim \exp\{\log\log N)^A\}$ with $A\ge2$. Then we see that
$$
L(N,\tau)\ge c\  N^{1-\s} \rho\left({\log N\over (\log\log N)^A}\right) \ge
c \ N^{1-\s} \exp\left({-c\, \log N\over (\log\log N)^{A-1}}\right)
$$
while the upper bound from Theorem 1.2 yields a better estimate
$$
\E\, \sup_{t \in \R} \big|\sum_{n\in {\cal E}_\tau}\e_n  n^{ -\s -it}\big|\le
C_\s \ N^{1/2-\s} \exp\left({- c\, \log N\over (\log\log N)^{A-1}}\right).
$$
The gap between the two upper bounds is polynomial.
One observes that $\ell_1$-estimate becomes even worse when $\tau$ decreases and approaches the critical zone.

\bigskip

{\bf Rudin-Shapiro polynomials.}
\
The upper bound in Theorem 1.2 is known to be related with Rudin-Shapiro
problem for Dirichlet polynomials. Let us recall first the classical
setting. For any trigonometric polynomial we have
$$
{\sum_{n=0}^{N-1}
|a_n| \over \sqrt{N}} \le
\sup_{t\in \R} | \sum_{n=0}^{N-1} a_n e^{int} |
\le 
\sum_{n=0}^{N-1} |a_n|.
\eqno(1.6)
$$
The arguments for getting the
lower bound are the inequality between the sup-norm and
$L_2$-norm, the
orthogonality of $(e^{int})_n$ and H\"older inequality. 

Rudin and Shapiro
constructed a fairly simple sequence $a_n\in\{-1,+1\}$ such that the
right
order of the lower bound is attained:
$$\sup_{t\in \R} |
\sum_{n=0}^{N-1} a_n e^{int} |
\le (2+\sqrt{2}) \sqrt{N+1}
\sim
(2+\sqrt{2}) \ {\sum_{n=0}^{N-1} |a_n| \over \sqrt{N}}.
$$
Consider
now the Dirichlet polynomials instead of the trigonometric ones. It is
known from
[KQ] and [Q3] that for any $(a_n)$
$$
\sup_{t\in \R} |
\sum_{n=0}^{N-1} a_n n^{it} | \ge 
\alpha_1 {\sum_{n=0}^{N-1} |a_n| \over
\sqrt{N}} 
\exp\{\beta_1 \sqrt{\log N\log\log N}\}.
$$
and for some
$(a_n)$
$$ 
\sup_{t\in \R} | \sum_{n=0}^{N-1} a_n n^{it} | \le 
\alpha_2
{\sum_{n=0}^{N-1} |a_n| \over \sqrt{N}} 
\exp\{\beta_2 \sqrt{\log N\log\log
N}\},
\eqno(1.7)
$$
with some universal constants
$\alpha_{1,2},\beta_{1,2}$.

Therefore, the lower bound for Dirichlet
polynomials is necessarily worse 
than in the classical case. Notice also
that the constuction
of example (1.7) in [Q3] is a probabilistic one; no
explicit example 
of Rudin-Shapiro type is known for Dirichlet
polynomials.
 It turns out that Theorem 1.2 generates a new family of
random polynomials
satisfying (1.7).

Indeed, take any $\s\in[0,{1\over
2})$ and choose $\tau$ 
in the optimal way. Namely, let
$$\log\tau\sim
({\log N\over 2})^{1/2} (\log\log N)^{1/2}.
$$
Set $a_n=\e_n n^{-\s}
1_{\{n\in {\cal E}_\tau \}}$.
It is easy to see that
$$
{ \sum_{n=0}^{N}
|a_n| \over \sqrt{N}}  =
{ \sum_{n\in {\cal E}_\tau}  n^{-\s} \over
\sqrt{N}}
\ge c\ N^{1/2-\sigma} \Psi_*(N,p_\tau),
$$
while by Theorem 1.2
we have the bound for the average 
of the left hand side in (1.7):
$$
\E
\sup_{t\in \R} | \sum_{n=0}^{N-1} a_n n^{it} | 
\le 
C_{\s}
N^{1/2-\sigma}\tau^{1/2}\Psi^*\Big({N \over
p_\tau^2},
p_{\tau}\Big)^{1/2}
$$
$$=
C_{\s}
N^{1/2-\sigma}
\exp\left\{
{1\over 2}
\left({\log N\over 2}\right)^{1/2}
(\log\log N)^{1/2}
+{1\over 2}
\log \Psi^*\Big({N \over
p_\tau^2},
p_{\tau}\Big) 
\right\}.
$$
Since by properties of Dickman
function,
$$\log \Psi^*\Big({N \over p_\tau^2},
p_{\tau}\Big) 
\sim
\log\rho\left(  {\log(N/p_\tau)\over \log p_\tau}\right)
\sim -\
{\log(N/p_\tau)\over \log p_\tau}
\log {\log(N/p_\tau)\over \log
p_\tau}
$$
$$
\sim 
-\  {\log N\over \log \tau}
\log {\log N\over \log
\tau}
\sim
- (2\log N)^{1/2} {(\log\log N)^{1/2}\over 2}=
- \left({\log
N\over 2}\right)^{1/2}\ (\log\log N)^{1/2}
$$
and by the same
arguments
$$\log \Psi^*\left(N, p_{\tau}\right) 
\sim 
- \left({\log N\over
2} \right)^{1/2} (\log\log N)^{1/2},
$$
we finally obtain
$$
\E
\sup_{t\in
\R} | \sum_{n=0}^{N-1} a_n n^{it} | 
\le 
{C_\s\over c} \
{ \sum_{n=0}^{N}
|a_n| \over \sqrt{N}} 
\exp\left\{\   \left({\log N\over 2} \right)^{1/2}
(\log\log N)^{1/2}
    \right\},
$$
as required in (1.7).
\bigskip

A particular case of this example with $\s=0$
was considered in [Q3]. Our calculation yields a slightly better
constant in the exponent. The question about the best possible constant
raised in [KQ] seems still to be open.
\bigskip
\bigskip

\noindent{\gum
2. Proof of the upper bound in Theorem 1.1.}
\medskip
\noindent The
principle of the proof of the upper bound
is as follows.
Once operated the
reduction to the study of a random
polynomial $Q$ on the multidimensional
torus by using
$(1.1)$, the proof then consists of two different steps
based on a
decomposition $Q=Q_1+ Q_2$. The study of the  supremum
of the
polynomial $Q_1$ is made by using the metric
entropy method.

The
investigation of the supremum of the polynomial
$Q_2$ is undertaken by
using first the contraction principle,
reducing the study to the one of a
complex valued Gaussian process.
The latter task is carried out by means of
Slepian's Comparison Lemma,
and by a careful study of the $L^2$-metric
induced by this process.

Now, we turn to the rigorous proof of the upper
bound and
introduce some notation.

We can represent ${\cal E}_\tau$ as the
union of disjoint sets
$$  E_j=\big\{ 2\le n\le N :  P^+(n)=p_j\big\},
\q
j=1,\ldots, \tau.
$$
\medskip

\noi
For $\uz \in \T^\tau$ we put
$$
Q(\uz)= \sum_{j=1}^\tau \sum_{n\in E_j}
\e_nn^{-\s}e^{2i\pi\langle
\ua(n),\uz\rangle}.
$$
By (1.1) we have
$$\sup_{t\in
\R}\big|\sum_{j=1}^\tau \sum_{n\in E_j} \e_n
n^{ -\s - it}\big| =\sup_{\uz
\in \T^\tau}\big|Q(\uz)\big|.
$$
Let $1\le \nu <\tau$ be fixed. Write
$Q=Q_1+Q_2$ where
$$Q_1(\uz)= \sum_{  P^+(n) \le p_{\nu }}
\e_nn^{-\s}
e^{2i\pi\langle\ua(n),\uz\rangle}, \qq
Q_2(\uz)= \sum_{p_{\nu
}<  P^+(n) \le  p_{\tau }}
\e_nn^{-\s}e^{2i\pi\langle
\ua(n),\uz\rangle}.
$$
First, evaluate the supremum of $Q_2$. Introduce the
following random process
$$   X^\e(\g) =\sum_{ \nu <j\le \tau}
\a_j
\sum_{n\in E_j}  \e_nn^{-\s} \b_{{n\over p_j }},
\qq \g\in
\G,
$$
where
$\g =\big((\a_j)_{\nu <j\le \tau}, (\b_m)_{1\le m\le N/2}\big)
$
and
$\G=\big\{ \g  :  |\a_j|\vee |\b_m|\le \!\!1,  \nu < \!j\le \tau,
1\!
\le m\le N/2\big\}$.
Writing  
$$
Q_2(\uz)= \sum_{ \nu <j\le \tau}
e^{2i\pi z_j  }\sum_{n\in E_j}
\e_nn^{-\s}e^{2i\pi\{\sum_{k\not = j}a_k(n)z_k+ [a_j(n)-1]z_j\}}
$$
and considering separately the imaginary and real parts of
$e^{2i\pi a_j(n)z_j  }$ and
$e^{2i\pi\sum_{k\not = j}a_k(n)z_k}$, easily
shows that
$$
Q_2(\uz)= X^\e(\g_1(\uz))+iX^\e(\g_2(\uz))+iX^\e(\g_3(\uz))+X^\e(\g_4(\uz)),
$$
where
$$
\g_1(\uz) = \big(\big( \cos(2\pi z_j)\big)_{\nu<j\le\tau}, 
                  \big(\cos(2\pi\sum_k a_k(m)z_k)\big)_{1\le m\le N/2} \big),
$$
$$
\g_2(\uz) = \big(\big( \sin(2\pi z_j)\big)_{\nu<j\le\tau}, 
                  \big(\cos(2\pi\sum_k a_k(m)z_k)\big)_{1\le m\le N/2} \big),
$$
etc. Therefore, we obtain
$$\sup_{\uz \in \T^\tau}\big|Q_2(\uz)\big|
\le 
4\sup_{\g \in\G}  \big|X^\e(\g)\big|.
$$
By the contraction principle ([K] p.16-17)

$$\E\, \sup_{\uz \in \T^\tau}\big|Q_2(\uz)\big|
\le  4 \ \sqrt{ \pi \over
2  }\
     \E\, \sup_{\g \in \G}\big|X(\g)\big|,
$$
where $\{X(\g), \g\in
\G\}$ is the same process as $X^\e(\g)$
except that the Rademacher random
variables $\e_n$ are replaced by
independent ${\cal N}(0,1)$ random
variables $\m_n$:
$$  X(\g) =\sum_{\nu<j\le \tau} \a_j \sum_{n\in E_j}
\m_nn^{-\s}
\b_{{n\over p_j}}.
$$
The problem now reduces to estimating the
supremum
of the real valued Gaussian process $X$. Towards  this aim,
we
examine the $L^2$-norm of its increments:
$$ \|X_\g-X_{\g'}\|_2^2   =
\sum_{ \nu <j\le \tau} \sum_{n\in E_j}
 n^{-2\s} \big[\a_j\b_{{n\over
p_j}}-\a'_j \b'_{{n\over
p_j}}\big]^2       \le 2\!\sum_{ \nu <j\le
\tau}\sum_{n\in E_j}
n^{-2\s} \big[(\a_j -\a'_j)^2+ (\b_{{n\over p_j}}-
\b'_{{n\over
p_j}})^2\big]  , $$
where we have used the identity
$\a_j\b_{{n\over p_j}}-\a'_j
\b'_{{n\over p_j}}= (\a_j -\a'_j)\b_{{n\over
p_j}}+
(\b_{{n\over p_j}}- \b'_{{n\over p_j}})\a'_j$.

\noi The "$\a$"
component part is easily controlled as follows,
$$\eqalign{\sum_{ \nu <j\le
\tau} \sum_{n\in E_j} n^{-2\s}
 (\a_j -\a'_j)^2
&\le
\sum_{\nu<j\le\tau}
(\a_j -\a'_j)^2p_j^{-2\s}
\sum_{m\le N/p_j} m^{-2\s} \cr
&\le
C_\s
\sum_{\nu<j\le\tau}
(\a_j -\a'_j)^2 \big({N^{1-2\s}\over p_j}\big)
.
\cr}
\eqno(2.1)$$
For the "$\b$" component part, we
have
$$\eqalign{
\sum_{\nu <j\le \tau} \sum_{n\in E_j}
{ (\b_{{n\over
p_j}}- \b'_{{n\over p_j}})^2\over n^{2\s}}
&\le  \sum_{m \le
N/p_\nu}
(\b_m-\b'_m)^2
 \big(\sum_{{\nu<j\le \tau\atop  mp_j\le N}}{1\over
(mp_j)^{2\s}}
 \big)
\cr
&:=
\sum_{m \le N/p_\nu} K_m^2
(\b_m-\b'_m)^2.
\cr}
\eqno(2.2)$$
Now we evaluate the coefficients $K_m$.
Consider two cases.

1)\ $m\le N/p_\tau$. Then $mp_j\le m p_\tau\le N$ for
all
$j\le \tau$ and, by using
the standard estimate  (see [HW], Theorem 8,
p.10)
$$ p_j \sim j \ \log j
\eqno(2.3)$$

\noi we
have
$$\eqalign{
K_m^2
&=
\sum_{\nu<j\le \tau } (mp_j)^{-2\s}
\le
m^{-2\s}
\sum_{j\le \tau } p_j^{-2\s}
\cr
&\le
C\ m^{-2\s} \sum_{j\le \tau } (j\,
\log j)^{-2\s}
=
C_\s \ m^{-2\s} \tau^{1-2\s} (\log
\tau)^{-2\s}
\cr
&\le
C_\s \ m^{-2\s} {\tau \over p_\tau^{2\s}}\
.
\cr}
$$
Thus
$$\eqalign{
\sum_{m\le N/p_\tau} K_m
&\le
C_\s {\tau^{1/2}
\over p_\tau^{\s}}\
\sum_{m\le N/p_\tau} \ m^{-\s}
\cr
&\le
C_\s
\left({N\over p_\tau}\right)^{1-\s}
{\tau^{1/2}\over
p_\tau^\s}
\cr
&=
{C_\s N^{1-\s} \tau^{1/2}\over p_\tau}\
\le {C_\s
N^{1-\s}\over \tau^{1/2}\log\tau}\ .
\cr}
$$

2)\ $N/p_\nu\ge m >
N/p_\tau$.
Then take a unique $k\in (\nu,\tau]$ such that
$N/p_k< m \le
N/p_{k-1}$. We have
$$\eqalign{
K_m^2 &= \sum_{\nu<j\le k-1 }
(mp_j)^{-2\s}
\le
m^{-2\s}\ \sum_{j\le k-1 } p_j^{-2\s}
\cr
&\le
C_\s\
m^{-2\s}\ \sum_{j\le k } (j\log j)^{-2\s}
\le
C_\s\ m^{-2\s}\  {k^{1-2\s}
\over (\log k)^{2\s}}
\cr
&\le C_\s m^{-2\s} \  {k \over p_k^{2\s}}
 \le
C_\s m^{-2\s} \  {k \over (N/m)^{2\s}}
\cr
&= C_\s \  {k \over N^{2\s}}\
.
\cr}
$$
Since $k\log k\le Cp_k\le C\ {N\over m}$, we have
$$ k\le C\
{N\over m} \ (\log({N\over m}))^{-1}.
$$
We arrive at
$K_m\le C_\s\
N^{-\s}   ({N\over m})^{1/2} \
(\log({N\over m}))^{-1/2}
$\ .
It follows
that
$$\eqalign{
\sum_{m\le N/p_\nu} K_m
&\le C_\s\  N^{-\s}\ \sum_{m\le
N/p_\nu}
({N\over m})^{1/2} \ (\log({N\over m}))^{-1/2}
\cr
&\le
C_\s\
N^{1-\s}\
\int_0^{1/p_\nu} u^{-1/2} \ (\log(1/u))^{-1/2} du
\cr
&\le
C_\s\
N^{1-\s}\
p_\nu^{-1/2}\ (\log p_\nu)^{-1/2}
\le
 {C_\s N^{1-\s}\over
\nu^{1/2} \log\nu}\ .
\cr}
$$

Now define a second Gaussian process by
putting for all $\g\in \G$
$$\eqalign{
Y(\g)& =
\sum_{\nu<j\le\tau}
\big({N^{1-2\s}\over p_j}\big)^{{1/ 2}}\a_j\xi'_j
+
\sum_{m\le N/p_\nu} K_m \ \b_m\xi''_m
 :=
 \ Y'_\g + Y''_\g ,
\cr}
$$
where $\xi'_i  $, $\xi''_j$  are independent
${\cal N}(0,1)$ random
variables. It follows from (2.1) and
(2.2) that for some suitable constant
$C_\s$, one has the
comparison relations: for all $\g, \g'\in
\G$,
$$\|X_\g-X_{\g'}\|_2\le C_\s \|Y_\g-Y_{\g'}\|_2.
$$
By virtue of the
Slepian
 comparison lemma (see [L], Theorem 4 p.190), since
$X_0=Y_0=0$,
we have
$$\E\, \sup_{\g\in \G}     |X_\g|\le
2 \E\,
\sup_{\g\in \G}      X_\g \le
2 C_\s \E\, \sup_{\g\in \G} Y_\g \le
2 C_\s
\E\, \sup_{\g\in \G} |Y_\g|.
$$
It remains to evaluate the supremum  of
$Y$. First of all,
$$\E\, \sup_{\g\in \G}|Y'(\g)|
\le N^{{1\over 2}-
\s}\sum_{\nu<j\le\tau}   p_j^{-1/2}
$$
By (2.3), we have
$$
\sum_{\nu<j\le\tau}   p_j^{-1/2}  \le
   \sum_{1<j\le\tau}   p_j^{-1/2}
\le
 {C \tau^{1/2} \over (\log\tau)^{1/2}}\ ,
$$
thus
$$
\E\, \sup_{\g\in
\G}|Y'(\g)|
\le C\ N^{{1\over 2}-\s}
\
{\tau^{1/2} \over (\log\tau)^{1/2}}
\ .
\eqno(2.4)$$

\noi To control the supremum of $Y''$, we use
our
estimates for the sums of $K_m$ and
write that
$$\eqalign{
 \E\,
\sup_{\g\in \G} |Y''(\g)|
&\le
   \sum_{m\le  N/p_\nu} K_m
\cr
&\le
C_\s
\left(
{ N^{1-\s} \over \nu^{1/2} \log\nu} +
{N^{1-\s}  \over \tau^{1/2}
\log\tau }
\right)
\le
{C_\s  N^{1-\s} \over \nu^{1/2} \log\nu}
\
.
\cr}
\eqno(2.5)
$$
\bigskip

Now, we turn to the supremum of $Q_1 $.
Towards this aim,
introduce the auxiliary Gaussian process
$$  \Upsilon
(\uz) =
\sum_{  P^+(n) \le p_\nu }   n^{-\s}
\big\{\t_n \cos 2\pi \langle
\ua(n),\uz\rangle +
\t_n'\sin 2\pi \langle \ua(n),\uz\rangle \big\}
,\qq
\uz\in \T^{\nu },
$$
where $\t_i$, $\t'_j$  are independent ${\cal N}(0,1)$
random variables. By symmetrization (see e.g. Lemma 2.3 p. 269 in
[PSW]),
$\displaystyle{\E\, \sup_{\uz \in \T^{\nu }}\big|Q_1(\uz)\big|\le
\sqrt{8\pi}
\E\, \sup_{\uz \in  \T^{\nu }}\big|\Upsilon (\uz)\big|}$, so
that we are again led to evaluating the supremum of a real valued
Gaussian process. For 
$\uz, \uz' \in \T^{\nu}$ put 
$\big\|\Upsilon(\uz)-\Upsilon(\uz)\big\|_2
 := \,d (\uz, \uz')$, 
and observe that
$$\eqalign{
\,d(\uz,\uz')^2  &= 4 \sum_{n: P^+(n)\le p_\nu}
  {1\over
n^{2\s}}
  \sin^2(\pi \langle\ua(n),\uz -\uz'\rangle)
\le
  4\pi^2\
\sum_{n: P^+(n) \le p_\nu}  {1\over n^{2\s}}
    |\langle\ua(n),\uz
-\uz'\rangle|^2
\cr
{} & \le
4\pi^2\
 \sum_{n: P^+(n)\le p_\nu}
n^{-2\s}
\left[\sum_{j=1}^\nu a_j(n) |z_j - z'_j|\right]^2
\cr
&  =
4\pi^2\ \sum_{n: P^+(n)\le p_\nu}   \sum_{j_1,j_2=1}^\nu
   a_{j_1}(n)
a_{j_2}(n) |z_{j_1} - z'_{j_1}|\ |z_{j_2} - z'_{j_2}|
n^{-2\s}
\cr
& =
4\pi^2\ \sum_{j_1,j_2=1}^\nu \sum_{n: P^+(n)\le p_\nu}
   a_{j_1}(n)
a_{j_2}(n)  |z_{j_1} - z'_{j_1}|\ |z_{j_2} - z'_{j_2}|
 n^{-2\s}
\cr
& \le
4\pi^2\
\sum_{j_1,j_2=1}^\nu  |z_{j_1} - z'_{j_1}|\ |z_{j_2} -
z'_{j_2}|
\sum_{b_1,b_2=1}^\infty b_1 b_2
\sum_{n\le N, a_{j_1}(n)=b_1,
a_{j_2}(n)=b_2}
 n^{-2\s}
\cr
& \le    4\pi^2\
\sum_{j_1,j_2=1}^\nu
|z_{j_1} - z'_{j_1}|\ |z_{j_2} - z'_{j_2}|
\sum_{b_1,b_2=1}^\infty
b_1 b_2
p_{j_1}^{-2 b_1\s} p_{j_2}^{-2 b_2\s}
\sum_{
{k\le N p_{j_1}^{-b_1}
p_{j_2}^{-b_2} \atop P^+(k)\le
p_\nu}}
k^{-2\s}
\cr
& \le
C_\sigma
N^{1-2\s}
\sum_{j_1,j_2=1}^\nu  |z_{j_1} - z'_{j_1}|\ |z_{j_2} -
z'_{j_2}|
\sum_{b_1,b_2=1}^\infty
b_1 b_2 p_{j_1}^{-2 b_1\s} p_{j_2}^{-2
b_2\s}
 [p_{j_1}^{-b_1} p_{j_2}^{-b_2}]^{1-2\s}
\cr
& =  C_\sigma
N^{1-2\s}
\sum_{j_1,j_2=1}^\nu  |z_{j_1} - z'_{j_1}|\ |z_{j_2} -
z'_{j_2}|
\sum_{b_1,b_2=1}^\infty
b_1 b_2 p_{j_1}^{-b_1}
p_{j_2}^{-b_2}
\cr
& =
C_\sigma N^{1-2\s}
\left\{
\sum_{j=1}^\nu  |z_{j} -
z'_{j}|
\sum_{b=1}^\infty
b \
p_{j}^{-b}
\right\}^2.
\cr}
\eqno(2.6)
$$
Thus,
$$ d(z,z')\le
C_\sigma
N^{1/2-\s}
\left\{
\sum_{j=1}^\nu  |z_{j} - z'_{j}|
\sum_{b=1}^\infty
b \
p_{j}^{-b}
\right\}.
\eqno(2.7)$$
\medskip

{\bf Remark.} In the middle of the long calculation,
we did not use the fact that the variable $k$
satisfies $P^+(k)\le p_\nu$. Actually,
this observation permits to introduce an extra factor related
do Dickman function, something like
$\rho(\log N/\log\nu)$.
This is helpful once $\nu$ is very small with respect to
$N$ (see the upper bound in Theorem 1.2).
\medskip

Now we explore the entropy properties of the
metric space $(\T^\nu,d)$. Towards
this aim, take $\e\in
(0,1)$ and cover $T^\nu$ by rectangular cells so that
if
$z$ and $z'$ belong to the same cell we
have
$$
|z_j-z'_j|\le
\cases{{\e\over \log\log\nu}&,\ $1\le
j\le\nu^{1/2}$,\cr
       \e                   &, $\nu^{1/2} <j\le
\nu$.
}
\eqno(2.8)$$
Thus, every cell is a product of two cubes of
different
size and dimension. The necessary number of cells $M(\e)$
is
bounded as follows,
$$ M(\e)\le \left({\log\log\nu\over
\e}\right)^{[\nu^{1/2}]}
\e^{-(\nu-[\nu^{1/2}])} = (1/\e)^\nu
(\log\log\nu)^{[\nu^{1/2}]}.
$$
Let us now evaluate the distance $d(z,z')$
for $z,z'$ satisfying
(2.8). By (2.7) we have
$$ d(z,z')\le
C_\sigma
N^{1/2-\s}
\left\{
d_1+d_2+d_3
\right\},
$$
where
$$ d_1= \sum_{j=1}^\nu
|z_{j} - z'_{j}|
\sum_{b=2}^\infty b \ p_{j}^{-b},
$$
$$ d_2
=\sum_{\nu^{1/2}<j\le \nu}  |z_{j} - z'_{j}| p_{j}^{-1},
$$
$$ d_3
=\sum_{j\le \nu^{1/2}}  |z_{j} - z'_{j}| p_{j}^{-1}.
$$
For any $j\ge 1$ we
have
$$  \sum_{b=2}^\infty b \ p_{j}^{-b} =
\sum_{b=2}^\infty b \ ({2\over
p_{j}})^{b} 2^{-b}
\le  ({2\over p_{j}})^{2} \sum_{b=2}^\infty b \
2^{-b}
= C  p_{j}^{-2}.
\eqno(2.9)
$$
Hence,
$$ d_1\le
\left(\sum_{j=1}^\nu
C  p_{j}^{-2}\right)\
\ \max_{j\le \nu}  |z_{j} - z'_{j}|
\le C
\e.
$$
Similarly,
$$ d_2\le
\left(\sum_{\nu^{1/2}<j\le\nu}
p_{j}^{-1}\right)\
\ \max_{\nu^{1/2}< j\le \nu}  |z_{j} - z'_{j}|
\le
C\
\left(\sum_{\nu^{1/2}<j\le\nu}  (j\log j)^{-1}\right)\
\e
$$
$$
\le C
\int_{\nu^{1/2}}^{\nu} {du\over u\log u}\ \e
= C (\log\log \nu-
\log({\log\nu\over 2})) \ \e
= C  (\log2) \, \e.
$$
Finally,
$$
d_3\le
\left(\sum_{j=1}^\nu p_{j}^{-1}\right)\
\ \max_{j\le \nu^{1/2}}
|z_{j} - z'_{j}|
\le
C \left(\sum_{j=1}^\nu (j\log j)^{-1}\right)\
{\e\over
\log\log\nu}
\le C\ \e.
$$
By summing up three estimates, we have
$d(z,z')\le C_\s N^{1/2-\s}
\e$ which enables the evaluation of the metric
entropy.

Let $\NNN\left(\T^\nu,d,u\right)$ be the minimal number of balls
of
radius $u$ that cover the space $(\T^\nu,d)$. We have

$$
\log\NNN\left(\T^\nu,d, C_\s N^{1/2-\s}\e\right) \le
\log M(\e) \le \nu
|\log\e|+ \nu^{1/2}\cdot \log\log\log\nu  .
$$
Observe also that
$$\|\Upsilon (\uz)\|_2\le C_\s N^{1/2-\s},\quad \uz\in
\T^{\nu}.
\eqno(2.10)
$$
Hence,
$D:=diam(\T^\nu,d)\le C_\s N^{{1\over
2}-\s}$, and by the
classical Dudley's  entropy theorem (see [L], Theorem 1
p.179),
for any fixed $\uz\in \T^\nu$
$$\eqalign{
\E\,  \sup_{\uz'\in
T^\nu} |\Upsilon (\uz')-\Upsilon(\uz)|&\le
C_\s\int_0^D
[\log\NNN(\T^\nu,d,u)]^{1/2}du \le  C_\s \int_0^{C_\s
N^{1/2-\s}}
[\log\NNN(\T^\nu,d,u)]^{1/2}du \cr &= C_\s N^{1/2-\s}
\int_0^1
[\log\NNN(\T^\nu,d,C_\s N^{1/2-\s}\e)]^{1/2}d\e \cr & \le
C_\s N^{1/2-\s}
\int_0^1\left[ \nu |\log\e|+ \log\log\log\nu \cdot
\nu^{1/2} \right]^{1/2}
d\e \cr &\le    C_\s N^{1/2-\s}\nu^{1/2}.
\cr}$$ Using again (2.10), we
have
$$
\E\,  \sup_{\uz'\in  T^\nu} |\Upsilon (\uz')|
\le    C_\s
N^{1/2-\s}\nu^{1/2}.
\eqno (2.11)
$$
The final stage of the proof provides
the optimal choice of
the parameter $\nu$ balancing the quantities (2.4),
(2.5),
and (2.11).
As suggests the Theorem's claim, we consider three
cases.

{\bf Case 1.}\ $N^{1/2}  \le
\tau\le N.$
Obviously, this case
contains the results of Hal\'asz and
Queff\'elec.
In this case we choose
$$\nu=
{\tau\over\log N}
$$
thus balancing (2.4) and (2.11). We obtain from both
terms
the bound
$   C_{\s} {N^{1/2-\sigma}\tau^{1/2}\over (\log
N)^{1/2}}
$
while the the term (2.5) is negligible.
The correctness
condition $\nu\le \tau$ is obvious.

{\bf Case 2.}\ $N^{1/2}(\log
N)^{-1}\le
\tau \le N^{1/2}.$
In this case we choose
$$\nu= N^{1/2} (\log
N)^{-1}
$$
thus balancing (2.5) and (2.11). We obtain from both terms
the
bound
$C_{\s} {N^{3/4-\sigma}
\over (\log N)^{1/2}}$
while thethe term
(2.4) is negligible.
The correctness condition $\nu\le \tau$ is obvious for
the range
under consideration.

{\bf Case 3.}\   $1\le \tau \le
N^{1/2}(\log N)^{-1}.$
Here we just set $\nu=\tau$. It means that we do not
need the
splitting of the polynomial in two parts. Formally,
the quantities
(2.4) and (2.5) are not necessary and we 
obtain the bound $C_{\s}
N^{1/2-\sigma}\tau^{1/2}$
directly from (2.11).

The upper bound is proved
completely.
\bigskip


\bigskip

\noindent{\gum 3. Proof of the lower bound in
Theorem 1.1.}
\medskip

\noi Let $\ud =\{d_n, n\ge 1\}$ be a   sequence of
reals. Recall
that by (1.1) we have
$$\sup_{t\in \R}\big|\sum_{j=1}^\tau
\sum_{n\in E_j} d_n\e_n
n^{ -\s - it}\big| =\sup_{\uz \in
\T^\tau}\big|Q(\uz)\big|.
$$
where
 $$ Q(\uz)= \sum_{j=1}^\tau \sum_{n\in
E_j}
   d_n \e_n n^{-\s} e^{2i\pi\langle \ua(n),\uz\rangle}.
$$

\noi
Consider the subset $\zz$ of $\T^\tau$ defined by
$$\zz=\Big\{ \uz=\{z_j,
1\le j\le \tau\} \quad : \quad \hbox{$z_j=0$,   
if $j\le \tau/2$,\quad
and
\quad $z_j\in\{0,1/2\}$, if $j\in(\tau/2,\tau]$} \Big\} .
$$

 Observe
that the imaginary part of $Q$ vanishes
on $\zz$, since for any $\uz\in
\zz$ and any $n$ it is true that

$$ e^{2i\pi\langle \ua(n),\uz\rangle} =
\cos(2\pi\langle
\ua(n),\uz\rangle) = (-1)^{2\langle
\ua(n),\uz\rangle}.
$$

Hence, $Q$ takes the following simple form on
$\zz$
$$
 Q(\uz)
= \sum_{\tau/2<j\le \tau}  \sum_{n\in E_j}
d_n \e_n n^{-\s}
{(-1)}^{2\langle \ua(n),\uz\rangle }
 . $$
 This is no longer a
trigonometric polynomial, but simply a finite
 rank Rademacher
process.

For $j\in(\tau/2,\tau]$ define

$$ {\cal L}_j=\Big\{n=p_j \,
\tilde n\ : \
  \tilde n\le {N\over p_{j}}\ \hbox{and}\
P^+(\tilde n)\le
p_{\tau/2}\Big\}. $$

Since
 $$  E_j \supset{\cal L}_j,
\q j=1,\ldots
\tau,
$$
the sets ${\cal L}_j$ are pairwise disjoint.

Put for $z\in
\zz$,
$$
 Q'(\uz)
= \sum_{  \tau/2<j\le \tau }\sum_{n\in  {\cal L}_j}
d_n \e_n
n^{-\s} {(-1)}^{2\langle \ua(n),\uz\rangle }
 .
 $$
We now recall a useful
fact.
\medskip

\noi {\gem Lemma 3.1.} {\it Let  $X=\{X_z, z\in Z\}$ and
$Y=\{Y_z,
z\in Z\}$ be two finite sets of random variables defined on
a
common probability space. We assume that $X$ and $Y$ are
independent and
that the random variables $Y_z$ are all centered.
Then}
$$ \E\sup_{z\in
Z}|X_z + Y_z| \ge  \E\sup_{z\in Z}|X_z  |.$$
{\it Proof.} Let $\F$ be the
$\sigma$-field generated by $Y$. Then
$$
\eqalign{ \E\sup_{z\in Z}|X_z +
Y_z| &= \E\left[\E\left(\sup_{z\in
Z} |X_z + Y_z|\ \big|\F\right)\right]
\cr & \ge \E\left[
\sup_{z\in Z}\big| \E ( X_z + Y_z \ \big|\F)
\big|\right] \cr & =
\E\left( \sup_{z\in Z}\big| X_z +\E Y_z \big|\right)
\cr &= \E
\sup_{z\in Z}\big| X_z \big|. \cr}
$$
\cqfd

Clearly, since
$\{Q(\uz)-Q'(\uz),\uz\in \zz\}$ and $\{
Q'(\uz),\uz\in \zz\}$ are
independent,
$$\E \sup_{\uz\in \zz} |Q(\uz)|
 \ge
  \E \sup_{\uz\in \zz}
\left|Q'(\uz) \right|.
$$

We now proceed to a direct evaluation of
$Q'(\uz)$ by proving
\medskip

\noi {\gem Proposition 3.2.} {\it There exists a universal
constant $c$ such that for any system of coefficients
$(d_n)$}
$$c\ \sum_{  \tau/2<j\le \tau }
\big|\sum_{n\in  {\cal L}_j}
d_n^2
\big|^{1/2} \le \E\, \sup_{\uz\in \zz}
\left|Q'(\uz)
\right|
\le
\sum_{  \tau/2<j\le \tau }
\big|\sum_{n\in  {\cal L}_j} d_n^2 \big|^{1/2}.
$$
\medskip

\noi {\it Proof.}
  For any $n\in
{\cal L}_j$, we have
$2\langle \ua(n),\uz\rangle= 2 z_j$, so that

$$\sum_{n\in  {\cal L}_j}
d_n\e_n  {(-1)}^{2\langle
\ua(n),\uz\rangle } ={(-1)}^{2z_j
}\sum_{n\in  {\cal L}_j} d_n\e_n(\o).  
$$
Thus
$$Q'(\uz)
= \sum_{  \tau/2<j\le \tau }{(-1)}^{2z_j
}\sum_{n\in  {\cal L}_j}
d_n\e_n(\o). 
$$ 
Let $\o\in \O$. We can select
$z_j=z_j(\o)=0$ or $1/2$,  $\tau/2<j\le \tau$, according to the sign
$+$ or $-$ of the sum $\sum_{n\in {\cal L}_j} d_n\e_n(\o)
n^{-\s}$. This implies that
$$\sup_{\uz\in \zz} \left|Q'(\uz) \right|
= \sum_{
\tau/2<j\le \tau }\big|\sum_{n\in  {\cal L}_j}
d_n\e_n  \big|.
$$
Consequently, by the Khintchine inequalities for Rademacher sums [KS]
$$
\eqalign{\E\, \sup_{\uz\in \zz} \left|Q'(\uz) \right| &=
\sum_{
\tau/2<j\le \tau }\E\, \big|\sum_{n\in  {\cal L}_j} d_n\e_n
\big|\ge c\  \sum_{  \tau/2<j\le \tau }\Big(\E\,
\big|\sum_{n\in  {\cal
L}_j} d_n\e_n  n^{-\s} \big|^2\Big)^{1/2}
\cr & =c \ \sum_{  \tau/2<j\le
\tau }\Big( \sum_{n\in  {\cal L}_j} d_n^2
\Big)^{1/2}.
\cr}
$$
The upper bound immediately follows from the Cauchy-Schwarz inequality.
\cqfd

\noi {\gem Corollary 3.3.} {\it If $(d_n$) is a multiplicative system, we have}
$$\E\, \sup_{\uz\in \zz}
\left|Q'(\uz) \right|\ge
c\ 
  \sum_{  \tau/2<j\le \tau }
d_{p_j}\Big( \sum_{\tilde n\le N/p_j\atop
  P^+(\tilde n)
  \le
p_{\tau /2}
} d^2_{\tilde n}\Big)^{1/2}
$$

 Now we can finish the proof of Theorem 1.1.
\medskip

\noi {\it Proof of the lower bound in Theorem 1.1.} If
$d_n\equiv n^{-\s}$, we get from the above  corollary
$$\eqalign{\E \sup_{\uz \in
\T^\tau}\big|\sum_{j=1}^\tau \sum_{n\in E_j}
  \e_n n^{-\s}e^{2i\pi\langle
\ua(n),\uz\rangle}\big|&\ge
  \E\, \sup_{\uz\in \zz} \left|Q'(\uz)
\right|\cr& \ge {C\over N^{ \s}}
\sum_{
\tau/2<j\le
\tau }  \#\big(  m\le
N/p_j : P^+(m)\le p_{\tau /2}  \big)^{1/2}\cr& = {C\over N^{
\s}}
\sum_{
\tau/2<j\le
\tau } \Psi({N\over p_j}   ,p_{\tau /2})^{1/2}
.
\cr}$$
Since
$$\eqalign{ \Psi({N\over p_j} ,p_{\tau /2})
&\ge
\Psi({N\over p_\tau} ,p_{\tau /2}) \cr & =
 {N\over p_\tau} \
\Psi^*({N\over p_\tau}\, ,p_{\tau /2})
\cr & \ge {c\ N\over \tau\log \tau }
\ \Psi^*({N\over p_\tau}\,
,p_{\tau /2}), \cr}
$$
we obtain
$$\eqalign{\E
\sup_{\uz \in \T^\tau}\big|\sum_{j=1}^\tau \sum_{n\in E_j}
\e_n n^{-\s} e^{2i\pi\langle \ua(n),\uz\rangle}\big| &\ge
{c\over N^{ \s}}\
{\tau\over 2} \ \left[ {c\, N\over \tau\log \tau
} \ \Psi^*({N\over
p_\tau}\, ,p_{\tau /2}) \right]^{1/2}\cr &= c \
N^{1/2-\s} \
\left(\tau\over\log \tau\right)^{1/2} \
\Psi^*({N\over p_\tau}\, ,p_{\tau
/2})^{1/2}, \cr} $$
as asserted.

\cqfd\medskip


\bigskip

\noindent{\gum 4. Proof of Theorem 1.2.} \

We need to prove the upper bound,
since the lower bound was obtained in
Theorem 1.1.
Moreover, we are only going to show how the calculations
concerning 
the upper bound of Theorem
1.1 should be corrected in order to
get an extra Dickman-type factor.

{\it Step 1. Some remarks on
semi-asyptotic formula for Dickman function.}

We discuss the so called
semi-asymptotic formula (see [BT])
$$\Psi(ax,y)=a^{\alpha(x,y)}
\Psi(x,y)
\left( 1+O(1/\bar u)\right)
\eqno(4.1)$$
where $\bar u=
\min\{\log x,y\}/\log y$ and
$$ \alpha(x,y)=
{\log(1+y/\log x)\over \log y}
=
1-{\log\log x\over \log y} +
{\log(1+\log x/y)\over \log
y}
$$
$$=
1-{\log\log x\over \log y} +
O\left({\log x\over y\log
y}\right).
$$

Since in our zone $y>\log x$, we have
$$O\left({\log
x\over y}\right)= O(1) = o(\log\log x).$$
Therefore $\alpha\le 1$ for $x$
large enough. We also see
that $\alpha\to 1$ when $x\to \infty$, hence
$\alpha\ge 2/3$
for all $x$ large enough.
We will use in the sequel that $
2/3\le \alpha\le 1.$
\bigskip

{\it Step 2. Main estimate and the adjustment of the previous proof.}
\medskip

We still use the notation
$\Psi^*(x,y)=x^{-1}$
but skip $y$ everywhere since $y=p_\nu$. In other words, we denote
$\Psi(x):=\Psi(x,p_\nu)$ and $\Psi^*(x):=\Psi^*(x,p_\nu)$.

Let $b^*=1$ for $b=1$ and $b^*=2b/3$ for $b=2,3,\dots$.
We will prove
that for all $b_1,b_2\ge 1$, $j_1,j_2\le\nu$
$$ \Psi\left({N\over
p_{j_1}^{b_1} p_{j_2}^{b_2}}\right)
\le C\ {N\over p_{j_1}^{b_1^*}
p_{j_2}^{b_2^*}}
\ \Psi^*\left({N\over
p_{\nu}^2}\right).
\eqno(4.2)$$

Once (4.2) is proved, the calculation from
(2.6) is updated as follows.
Let denote $D_j=|z_{j} - z'_{j}|$. Then
$$
d(z,z')^2\le C\ \sum_{j_1,j_2\le \nu} D_{j_1}
D_{j_2}
\sum_{b_1,b_2=1}^{\infty}b_1 b_2
p_{j_1}^{-2b_1\sigma}
p_{j_2}^{-2b_2\sigma} \Psi\left({N\over
p_{j_1}^{b_1} p_{j_2}^{b_2}}\right)
\ \left({N\over 
p_{j_1}^{b_1}
p_{j_2}^{b_2}}\right)^{-2\sigma}
$$
$$=C\ N^{-2\sigma} \ \sum_{j_1,j_2\le
\nu} D_{j_1} D_{j_2}
\sum_{b_1,b_2=1}^{\infty}b_1 b_2 
 \Psi\left({N\over
p_{j_1}^{b_1} p_{j_2}^{b_2}}\right)
$$
$$\le C\ N^{1-2\sigma} \ 
\Psi^*
\left({N\over p_{\nu}^{2}}\right)
\sum_{j_1,j_2\le \nu} D_{j_1}
D_{j_2}
\sum_{b_1,b_2=1}^{\infty}b_1 b_2
p_{j_1}^{-b_1^*}
p_{j_2}^{-b_2^*}
$$
$$= C\ N^{1-2\sigma} \ 
\Psi^*
\left({N\over p_{\nu}^{2}}\right)
\left\{\sum_{j\le \nu}
D_{j}
\sum_{b=1}^{\infty}b p_{j}^{-b^*}
\right\}^2.
$$
Now everything
continues as in the proof of Theorem 1.1 but with an extra factor
$\Psi^*
\left({N\over p_{\nu}^{2}}\right)$. The minor change 
(corresponding to
(2.9) is that
$$\sum_{b=2}^{\infty} b p_{j}^{-b^*}
=\sum_{b=2}^{\infty}
b\left({2\over p_{j}}\right)^{b^*} 2^{-b^*}\le
\left({2\over
p_{j}}\right)^{4/3}   \sum_{b=2}^{\infty} b 2^{-b^*}
={C\over
p_{j}^{4/3}},
$$
hence still
$$ d_1\le \sum_{j=1}^{\nu} {C\over
p_{j}^{4/3}} \ 
\max_j D_j \le C\epsilon
$$ 
\bigskip

{\it Step 3. The proof of inequality $(4.2)$}.

We consider three cases

1. $b_1,b_2\ge 2$. \
By applying (4.1) with $x={N\over  p_{j_1}^{b_1} 
p_{j_2}^{b_2}}$ and
$a=p_{j_1}^{b_1} p_{j_2}^{b_2}$, we get
$$ \Psi(N)= \left(p_{j_1}^{b_1}
p_{j_2}^{b_2}\right)^\alpha
\Psi\left({N\over p_{j_1}^{b_1}
p_{j_2}^{b_2}}\right)
\left(1+O(1/\bar u)\right).
$$
Once $\bar u$ is large
enough and $\alpha\ge 2/3$ we have
$$
\Psi\left({N\over p_{j_1}^{b_1}
p_{j_2}^{b_2}}\right)\le
C\ \Psi(N) \left(p_{j_1}^{b_1}
p_{j_2}^{b_2}\right)^{-2/3}.
$$
Similarly, we pass from $\Psi(N)$ to
$\Psi\left({N\over p_{\nu}^{2}}\right)$.
By using $\alpha\le 1$, we
have
$$
\Psi(N)=[p_{\nu}^{2}]^\alpha \Psi\left({N\over
p_{\nu}^{2}}\right)
\left(1+O(1/\bar u)\right)
\le C p_{\nu}^{2}\
\Psi\left({N\over p_{\nu}^{2}}\right)
=C N \Psi^* \left({N\over
p_{\nu}^{2}}\right).
$$
By combining two estimates we
get
$$\Psi\left({N\over p_{j_1}^{b_1} p_{j_2}^{b_2}} \right)\le
C   \
\left(p_{j_1}^{b_1} p_{j_2}^{b_2} \right)^{-2/3}
N \Psi^* \left({N\over
p_{\nu}^{2}}\right),$$
as required.
\medskip

2. $b_1=b_2=1$.
By applying
(4.1) with $x={N\over  p_{j_1} p_{j_2}}$ and
$a={p_{\nu}^{2}\over  p_{j_1}
p_{j_2} }$, we get, using
$\alpha\le 1$,
$$
\Psi \left( { N\over p_{j_1}
p_{j_2} } \right)
\le
C \left({ p_{\nu}^{2} \over p_{j_1} p_{j_2}}
\right)^{\alpha} 
\Psi \left({N\over p_{\nu}^{2}}\right)
\le
C\ {
p_{\nu}^{2} \over p_{j_1} p_{j_2}}  \
\Psi \left({N\over
p_{\nu}^{2}}\right)
= C\  {N \over p_{j_1} p_{j_2}}  \
\Psi^* \left({N\over
p_{\nu}^{2}}\right),
\eqno(4.3)
$$
as required.
\medskip

3. $b_1=1,b_2\ge
2$.
By applying (4.1) with $x={N\over  p_{j_1} p_{j_2}^{b_2}}$
and
$a={p_{j_2}^{b_2} }$, we get, using
$\alpha\ge 2/3$,
$$ 
\Psi({N\over
p_{j_1}})= \left(p_{j_2}^{b_2}\right)^\alpha
\Psi\left({N\over p_{j_1}
p_{j_2}^{b_2}}\right)
\left(1+O(1/\bar
u)\right),
$$
hence
$$
\Psi\left({N\over p_{j_1} p_{j_2}^{b_2}}\right)
\le
C p_{j_2}^{-2b_2/3}
\Psi\left({N\over p_{j_1}}\right)
=
C
p_{j_2}^{-b_2^*}
\Psi\left({N\over p_{j_1}}\right).
$$
Yet, letting $p_{j_2}=1$ in (4.3), we have
$$
\Psi\left({N\over p_{j_1}}\right)\le
C\ {N \over p_{j_1}}  \
\Psi^* \left({N\over p_{\nu}^{2}}\right),
$$
and we are done with case 3. Therefore, the proof of (4.2) is complete.
\bigskip

We finish the section by giving a proof of $(1.4')$. Only the upper bound needs 
a proof. Fix some large integer $M$. Let $\{g_n, n\ge 1\}$ be a sequence of
independent ${\cal N}(0,1)$ distributed random variables. 
By contraction principle, there is an absolute constant $C$ such that 
$$   \E\, \sup_{N\le M}
\sup_{t \in \R}  {|\sum_{n=2}^N \e_n  n^{-\s  - it}| \over  N^{1-\s}     
(\log  N)^{-1}   }  \le    C  \  \E\, \sup_{N\le M}
\sup_{t \in \R}  {|\sum_{n=2}^N g_n  n^{-\s  - it}| \over  N^{1-\s}     
(\log  N)^{-1}   }.
$$
We now need the following inequality (see [W1] p.451) which is a simple
consequence of Borell-Sudakov-Tsirelson inequality: 
{\it if $G_1,\ldots, G_N$ are Gaussian random vectors 
with values in a separable Banach space $(B,\|\cdot\|)$, then 
$$\E \ \sup_{1\le k\le N} ||G_k||\le C 
\left\{ \sup_{1\le k\le N  } \E ||G_k|| + \E \sup_{1\le k\le N} \s_k
|\zeta_k|\right\}
$$ 
where $\s_k 
=\sup_{f\in B^*,\ ||f||\le 1} \big(\E \langle f,G_k \rangle ^2\big)^{1/2}$, 
$k=1,\ldots, N$,
$\{\zeta_k, 1\le  k\le N\} $ is a sequence of independent 
${\cal N}(0,1)$ distributed random variables, and $C$ is a universal
constant.}
\smallskip

Applying this inequality gives
$$\eqalign{ \E\, \sup_{N\le M}
\sup_{t \in \R}  {|\sum_{n=2}^N g_n  n^{-\s  - it}| \over  N^{1-\s}     
(\log  N)^{-1}   }&\le C \sup_{N\le M}\E\,  
\sup_{t \in \R}  {|\sum_{n=2}^N g_n  n^{-\s  - it}| \over  N^{1-\s}     
(\log  N)^{-1}   } + C \E\, \sup_{N\le M}
    | \zeta_N\s_N|  \cr & \le C_\s +   C \E\, \sup_{N\le M}
    | \zeta_N\s_N|  ,\cr}
$$
where 
$$\s_N\le C{\sup_{t \in \R}\|\sum_{n=2}^N g_n  n^{-\s  - it}\|_{2}\over  N^{1-\s}     
(\log  N)^{-1}   } \le C{\big(\sum_{n=2}^N    n^{-2\s 
}\big)^{1/2}\over  N^{1-\s}     (\log  N)^{-1}   } \le 
C_\s {N^{1/2- \s}\over  N^{1-\s}     (\log  N)^{-1}   }= C_\s { \log N\over  N^{1/2 }    
   }.
$$
It is an obvious fact  that 
$\E\, \sup_{N\le M} {| \zeta_N |\log N \over  N^{1/2} }$ 
is bounded uniformly in $M$ 
by some absolute constant. 
So that, there exists a constant $C_\s$ depending on
$\s$ only, such that for any $M$ 
$$\E\, \sup_{N\le M}
\sup_{t \in \R}  {|\sum_{n=2}^N g_n  n^{-\s  - it}| \over  N^{1-\s} 
(\log  N)^{-1}   }\le C_\s. 
$$
The claimed result follows immediately. 

Note to conclude that the same argument applies to our upper bounds 
results with minor modifications (by introducing suitable blocks). 
\bigskip\bigskip

\noindent{\gum 5. Other results.}
\medskip

\noi In this section we test our technique on some other sets of coefficients.

Let $\{d_n, n\ge 1\}$ be a
sequence of multiplicative weights: $d_{nm}= d_nd_m$ whenever $n,m$
are coprimes. Denote
$$B_m= \sum_{2\le n\le m} d_n^2.
\eqno(5.1)
$$  
By choosing $\tau=\mu := \pi(N) $ in the lower bound of Proposition 3.2, we get
$$\eqalign{\E \sup_{\uz \in \T^\m}\big|\sum_{n =2}^N
d_n \e_nn^{-\s}e^{2i\pi\langle \ua(n),\uz\rangle}\big|&\ge \E\,
\sup_{\uz\in \zz} \left|Q''(\uz) \right| \ge CN^{-\s}  \sum_{
\m/2<j\le \m}  d_{p_j}\Big( \sum_{\tilde n\le N/p_j\atop P^+(\tilde n)\le p_{\m /2} }
   d^2_{\tilde n}\Big)^{1/2}.
\cr}$$
Note that for large $N$ in the case $\tau=\mu$ the sets $ {\cal L}_j$ reduce to  
$ \big\{n=p_j \, \tilde n\ : \  \tilde n\le {N\over p_{j}}
 \big\}.$
Indeed,
  if $\tilde n\le {N\over p_{j}}$ and if there is an $s\ge \m/2$
  such that $p_s|\tilde n$, then this implies that 
$$N\ge p_j p_s \ge
 p_{\m/2}^2 \sim( \m\log \m)^2/4\sim N^2/4,
 $$
 which is impossible for large $N$. Thus necessarily $P^+(\tilde n)\le p_{\m/2}$. 
Thereby,
$$\eqalign{\E \sup_{\uz \in \T^\m}\big|\sum_{n =2}^N
d_n \e_nn^{-\s}e^{2i\pi\langle \ua(n),\uz\rangle}\big|&  \ge CN^{-\s}  \sum_{
\m/2<j\le \m}  d_{p_j}\Big( \sum_{\tilde n\le N/p_j\  }
   d^2_{\tilde n}\Big)^{1/2}
\cr&=CN^{-\s}  \sum_{
\m/2<j\le \m}  d_{p_j}B_{N/p_j}^{1/2} .
\cr}$$
We have obtained
\medskip

\noi {\gem Proposition 3.4.} {\it There exists a universal constant
$C,N_0$ such that for any $0\le \s < 1/2$,   any integer $N\ge N_0$ and any
multiplicative sequence of weights $(d_n)$
$$\E
\sup_{t\in \R} \big|\sum_{n=2}^N\e_n \aa_n n^{-\s - it}\big|\ge C
N^{-\s}  \sum_{ \m/2<j\le \m}  d_{p_j}B_{N/p_j}^{1/2}, 
$$ 
where $B_m$ is defined in $(5.1)$.}
\medskip

Apply this to the case $d_n=d(n)$, where $d(n)= \sharp\{ d: d|n\}$
is the divisor function.  Although  these weights are very
irregular,  their sums behave regularly, in particular,
$$
 \sum_{n=1}^N d^2(n) \sim  ({N\over \pi^2}) \log^3 N .
$$
as $N$ tends to infinity. The last estimate immediately provides
$ B_m \sim  ({m/ \pi^2}) \log^3 m$, hence
 (noticing that $d_{p_j}=2$ and $\m\sim N/\log N)$)
$$ \sum_{
\m/2<j\le \m}  d_{p_j}B_{N/p_j}^{1/2} 
\sim  
\sum_{ \m/2<j\le \m}    ({2 N/p_j  \pi^2})^{1/2} \log^{3/2} {N\over p_j } 
=
{2N^{1/2}\over \pi} \sum_{
\m/2<j\le \m}  {1 \over p_j^{1/2} }    \log^{3/2} {N\over p_j }  $$

$$\sim 
{2N^{1/2}\over \pi} \sum_{ \m/2<j\le \m}  { \log^{3/2} 
              {N\over j\log j }\over (j\log j)^{1/2} }
\approx 
N^{1/2} \sum_{
\m/2<j\le \m}  { 1\over (j\log j)^{1/2} }
\approx
N^{1/2}
 { \m^{1/2}\over ( \log \m)^{1/2} }   
\sim  
{N\over   \log N} \ . 
$$

Now, let $\{P_k, k\in K \}$ be a finite set of mutually coprime numbers.
Consider the set of integers
$$E= \big\{n: \ n=\prod_{k\in K} P_k^{\a_k}, \qq \a_k\in\{0,1\}\big\}  $$
and the associated Dirichlet polynomial
$$ D_E(t)=\sum_{n\in E}\e_n   n^{-\s - it}=\sum_{n=2}^N\e_n \chi_E(n)  n^{-\s - it} , $$
where  $N=\prod_{k\in K} P_k$.
We prove the following.
\medskip

\noi {\gem Proposition 3.5.}\ {\it There exists a universal constant
$C$ such that, for any $\s\ge 0$ and any $\{P_k, k\in K \}$
$$\E\sup_{t\in \R} | D_E(t) |\ge  
C\prod_{k\in K} \Big(1+  P_k^{-2\s }  \Big)^{1/2}
\sup_{G\subseteq K }
 {\sum_{j\in G}P_{j }^{- \s }\over
 \displaystyle{ \prod_{ k\in G}  \Big(1+  P_k^{-2\s }  \Big)^{1/2}} }
  .$$
}

\noi{\it Proof.} By (1.1) we have
$$\sup_{t\in \R}\big| D_E(t)\big| =\sup_{\uz \in \T^\m}\big|Q(\uz)\big|.
$$
where $\mu=|K|$ and
 $$ Q(\uz)= \sum_{n=2}^N  \chi_E(n)
    \e_n n^{-\s}e^{2i\pi\langle \ua(n),\uz\rangle}.
$$

Let $A\subset K $ and $B= K \backslash A$. We assume that both $A$
and $B$ are nonempty sets.   Define for $j\in B$,
$${\cal B}_j=\{n\in E :
\a_k=0\ {\rm if}\ k\in B, k\not =j  , \
\a_j=1
\}$$ 
and $\zz\subset \T^\m$ by
$$\zz=\Big\{ \uz=\{z_k,  1\le k\le 2r\} \quad : \quad \hbox{$z_k=0$,
if $k\in A$,\quad and
$z_k\in \{0,1/2\}$ if $k\in B$\Big\}}.$$
For $j\in B$,  $n\in {\cal B}_j$  and $z\in \zz$, we have
$2\langle \ua(n),\uz\rangle= 2\sum_{k\in K} \a_k z_k= 2 z_j=\pm 1$, so that
similarly to our previous lower bound
$$\sup_{\uz \in \zz}\big|Q(\uz)\big|\ge
\sum_{j\in B}\big|\sum_{n\in {\cal B}_j}\e_nn^{-\s}\big|, $$
almost surely. Hence
$$\eqalign{
\E \sup_{\uz \in \zz}\big|Q(\uz)\big|
&\ge 
C \ \sum_{j\in B} \Big( \E\big|\sum_{n\in {\cal B}_j}\e_n
n^{-2\s}\big|^2 \Big)^{1/2}
= 
C\sum_{j\in B} 
P_{j }^{-\s }
 \Big( \sum_{(\a_k)_{k\in A}\in\{0,1 \}^A   }
 \prod_{ k\in A}P_k^{-2\s\a_k} \Big)^{1/2} \cr 
&=
 C\prod_{ k\in A}\Big(1+  P_k^{-2\s }
\Big)^{1/2} \Big\{  \sum_{j\in B}  P_{j }^{- \s } \Big\}.
\cr}
$$
Therefore 
$$\eqalign{
\E\sup_{t\in \R}\big| D_E(t)\big|
&\ge 
C\ \sup_{A\subseteq K, A\not=K}  \prod_{ k\in A}
\Big(1+  P_k^{-2\s }
\Big)^{1/2}
\Big\{\sum_{j\in A^c}  P_{j }^{- \s }\Big\} 
\cr 
&=
C\prod_{k\in K} \Big(1+  P_k^{-2\s }  \Big)^{1/2}
\sup_{A\subseteq K, A\not=K}
 {\sum_{j\in A^c}P_{j }^{- \s }\over
 \displaystyle{ \prod_{ k\in A^c}  \Big(1+  P_k^{-2\s }  \Big)^{1/2}} }.
 \cr 
}$$
\cqfd
\medskip
{\bf Acknowledgements.}\ The work of the first mentioned author was 
supported by grants RFBR 05-01-00911 and INTAS 03-51-5018. 
He is also grateful for hospitality to the L.Pasteur University (Strasburg)
where this research has been done.
\medskip

\noi {\gum References}
\smallskip

\noi [{B}]  {Bohr H.} [1952]  {\sl Collected Mathematical
Works}, Copenhagen.

\noi [BH] {Bohnenblust H.F., Hille E.} [1931] {\sl On the
absolute convergence of Dirichlet  series},  Ann. Math.
{\bf 2(32)}, 600--622.

\noi [BKQ] Bayart F., Konyagin S. V., Queff\'elec H. [2003/2004] 
{\sl Convergence almost everywhere and divergence almost everywhere of Taylor and Dirichlet
series} Real Anal. Exchange  {\bf 29}, no. 2, 557--586. 

\noi [BT] de la Bret\`eche, R.; Tenenbaum, G. [2005] {\sl Entiers friables: 
in\'egalit\'e de Tur\'an-Kubilius et applications.} 
Invent. Math. {\bf 159}, no. 3, 531--588. 
 
\noi [C]  Clarke  L. E. [1969] {\sl Dirichlet series with independent and 
identically disturbed coefficients}, Proc. Cambridge Philos. Soc. {\bf 66} 393--397.

\noi [DC] Dvoretzky A.,  Chojnacki  H. [1947] {\sl Sur les changements de signe d'une 
s\'erie \`a termes complexes}. C. R. Acad. {\bf 222}, 515--518.

\noi [DE1] Dvoretzky A.,  Erd\" os  P. [1955] {\sl On power series diverging 
everywhere on the circle of convergence}. Michigan Math. J. {\bf 3}, 31--35.
 
\noi [DE2] Dvoretzky A.,  Erd\" os  P. [1959] {\sl Divergence of random power 
series}, Michigan Math. J. {\bf 6},  343--347.  
 
\noi [H] {Helson H.} [1967] {\sl  Foundations of the theory of 
Dirichlet series}, Acta Math., {\bf 118}, 61--77.
 
\noi [Ha1] {Hal\'asz G.},   private communication to H.Queff\'elec, see [Q1]. 
 
\noi [Ha2] {Hal\'asz G.}   [1983]   {\sl  On random multiplicative functions},  Hubert Delange 
colloquium (Orsay 1982), Pub. Math. Orsay, {\bf 83-4}, Univ. Paris XI, Orsay,   74-96. 

\noi [Har] { Hartman  P.} [1939] {\sl On Dirichlet series involving random coefficients},  
Amer. J. Math. {\bf 61}, 955--964.

\noi [HR] {Hardy G.H., Riesz M.} [1915] {\sl The general theory of Dirichlet's series}, 
Cambridge Tracts in Math. and Math. Phys. {\bf 18}.

\noi [HS] {Hedenmalm  H., Saksman  E.} [2003] {\sl  Carleson's convergence theorem 
for Dirichlet series} Pacific J. Math. {\bf 208}  no. 1, 85--109.

\noi  [{HW}] {Hardy G.H., Wright E.M.} [1979] {\sl An Introduction 
to the Theory of Numbers}, Oxford University Press, Clarendon Press, Fifth ed.

\noi  [{K}] {Kahane J. P.} [1968] {\sl Some random series of functions}, 
D. C. Heath and Co. Raytheon Education Co., Lexington, Mass.
 
\noi [KQ] {Konyagin S.V., Queff\'elec H.} [2001/2002] {\sl
The translation ${1\over 2}$ in the theory of Dirichlet series},
Real Anal. Exchange {\bf 27}(1), 155--176.

\noi [KS] {Kashin B.S., Saakyan A.A.} [1989]:   {\sl Orthogonal Series}, 
Translations of Mathematical Monographs {\bf 75}, American Math. Soc. 

\noi [L] {Lifshits M.A.} [1995] {\sl Gaussian Random Functions},
Kluwer, Dordrecht.

\noi [PSW] {Peskir G., Schneider D., Weber M.} [1996] {\sl
Randomly weighted series of contractions in Hilbert spaces}, 
Math. Scand. {\bf 79}, 263--282.

\noi [Q1] {Queff\'elec H.} [1980] {\sl Propri\'et\'es presque
s\^ures et quasi-s\^ures des s\'eries de Dirichlet et des
produits d'Euler},  Can. J. Math.  {\bf XXXII} no. 3, 531--558.

\noi [Q2] {Queff\'elec H.} [1983] {\sl Sur une estimation
probabiliste li\'ee \`a l'in\'egalit\'e de Bohr}, In: 
Harmonic analysis: study group on translation-invariant 
Banach spaces,  Exp. No. 6, 21 pp., Publ. Math. Orsay, 84-1,
Univ. Paris XI, Orsay, 1984. 

\noi [Q3] {Queff\'elec H.} [1995] {\sl H. Bohr's vision of
ordinary Dirichlet series; old and new results},  J. Analysis  {\bf 3}, 
p.43-60.

\noi [STY] {Sun D.,  Tian F., Yu J.R.} [1998] {\sl Sur les s\'eries
al\'eatoires de  Dirichlet},  C. R. Acad. Sci. Paris S\'er. 1 {\bf
326}, 427--431.

\noi [T] {Tenenbaum G.} [1990] {\sl Introduction  \`a
la  th\'eorie analytique et probabiliste des nombres},
Revue de l'Institut Elie Cartan  {\bf 13}, D\'epartement 
de Math\'ematiques de l'Universit\'e de Nancy I.

\noi [W1]\ {Weber M.} [2000]  {\sl  Estimating random
polynomials by means of metric entropy methods},
Math. Inequal. Appl. {\bf 3},  no. 3, 443--457.

\noi [W2]\ {Weber M.} [2006]  {\sl  On a stronger form of
Salem-Zygmund inequality for random polynomials},  
Periodica Math. Hung.  {\bf 52}, No. 2, 73--104.

\noi [Y1] {Yu J.R.} [1978] {\sl Some properties of random Dirichlet series},  
Acta Math. Sinica {\bf 21}, 97--118.

\noi [Y2] {Yu J.R.} [1985] {\sl Sur quelques s\'eries gaussiennes de  Dirichlet},  
C. R. Acad. Sci. Paris S\'er. 1 {\bf 300}, 521--522.

\noi [Y3] {Yu J.R.} [1995] {\sl Dirichlet spaces and random
Dirichlet series},  J. Analysis  {\bf 3}, 61--71.
\smallskip

\noi {\phh Mikhail   Lifshits,  St.Petersburg State University,
Department of   Mathemmatics and Mechanics, 198504, Bibliotechnaya pl, 2, 
Stary Peterhof, Russia. E-mail:\ \tt lifts@mail.rcom.ru}
  
\noi {\phh Michel  Weber, \noi  Math\'ematique (IRMA),
Universit\'e Louis-Pasteur et C.N.R.S.,   7  rue Ren\'e Descartes, 
67084 Strasbourg Cedex, France. E-mail: \  \tt weber@math.u-strasbg.fr}

\bye